\documentclass[a4paper,11pt,twoside,english]{article}
\usepackage{babel}
\usepackage{latexsym,amsfonts}
\usepackage{amsthm}
\usepackage{amsmath}
\usepackage{mathtools}
\usepackage{paralist}

\usepackage{lastpage}
\usepackage{fancyhdr}
\fancypagestyle{plain}{

  \fancyhead[L,C]{}
  \fancyhead[R]{\today}
  \fancyfoot[R,L]{}
  \fancyfoot[C]{\thepage \ of \pageref{LastPage}}
  }

\title{Ball generated property of direct sums of Banach spaces}
\author{Jan-David Hardtke}
\date{}

\providecommand{\sm}{\setminus}
\providecommand{\ssq}{\subseteq}
\providecommand{\N}{\ensuremath{\mathbb{N}}}

\providecommand{\R}{\ensuremath{\mathbb{R}}}
\providecommand{\C}{\ensuremath{\mathbb{C}}}
\providecommand{\A}{\ensuremath{\mathcal{A}}}
\providecommand{\B}{\ensuremath{\mathcal{B}}}
\providecommand{\eps}{\ensuremath{\varepsilon}}

\providecommand{\keywords}[1]{
  {\let\thefootnote=\relax
  \footnote{{\em Keywords}: #1}}
  \addtocounter{footnote}{-1}
  }

\providecommand{\AMS}[1]{
  {\let\thefootnote=\relax
  \footnote{{\em AMS Subject Classification} (2010): #1}}
  \addtocounter{footnote}{-1}
  }

\providecommand{\address}{
  {\sc \noindent Department of Mathematics \\
  Freie Universit\"at Berlin \\
  Arnimallee 6, 14195 Berlin \\
  Germany \\}
  }

\DeclarePairedDelimiter{\set}{\lbrace}{\rbrace}

\DeclarePairedDelimiter{\abs}{\lvert}{\rvert}
\DeclarePairedDelimiter{\norm}{\lVert}{\rVert}

\theoremstyle{definition}
\newtheorem{definition}{Definition}[section]
\newtheorem*{definition*}{Definition}

\newtheorem*{example*}{Example}

\newtheorem*{remark*}{Remark}
\theoremstyle{plain}
\newtheorem{lemma}[definition]{Lemma}
\newtheorem*{lemma*}{Lemma}
\newtheorem{proposition}[definition]{Proposition}
\newtheorem*{proposition*}{Proposition}
\newtheorem{theorem}[definition]{Theorem}
\newtheorem*{theorem*}{Theorem}
\newtheorem{corollary}[definition]{Corollary}
\newtheorem*{corolary*}{Corollary}

\newenvironment{Proof}[1][\proofname]{\begin{proof}[#1] \setlength{\parindent}{0pt}}{\end{proof}}
\newenvironment{Abstract}{\centering\begin{minipage}{0.8\textwidth} \noindent \small {\sc Abstract.}}{\end{minipage}\par}

\usepackage{color}
\definecolor{darkgreen}{rgb}{0,0.5,0}

\numberwithin{equation}{section}
\newtagform{colored}[\color{blue}]{\color{blue}(}{\color{blue})}
\usetagform{colored}

\usepackage[colorlinks,linkcolor=blue,citecolor=red,urlcolor=darkgreen]{hyperref}
\providecommand{\email}{{\it E-mail address:} \href{mailto:hardtke@math.fu-berlin.de}{\tt hardtke@math.fu-berlin.de}}

\usepackage{amsrefs}

\begin{document}

\maketitle

\begin{Abstract}
\noindent A Banach space $X$ is said to have the ball generated property (BGP) if every closed, bounded, convex 
subset of $X$ can be written as an intersection of finite unions of closed balls.\par
In \cite{basu} S. Basu proved that the BGP is stable under (infinite) $c_0$- and $\ell^p$-sums for $1<p<\infty$.
We will show here that for any absolute, normalised norm $\norm{\cdot}_E$ on $\R^2$ satisfying a certain smoothness
condition the direct sum $X\oplus_E Y$ of two Banach spaces $X$ and $Y$ with respect to $\norm{\cdot}_E$ enjoys the BGP
whenever $X$ and $Y$ have the BGP.
\end{Abstract}
\AMS{46B20}

\section{Introduction}\label{sec:intro}
Let $X$ be a real Banach space. For $x\in X$ and $r>0$ we denote by $B_r(x)$ the closed
ball with center $x$ and radius $r$. The closed unit ball $B_1(0)$ is simply denoted by $B_X$, while 
$S_X$ stands for the unit shpere. Finally, $X^*$ denotes the dual space of $X$.\par
$X$ is said to have the ball generated property (BGP) if every closed, bounded, convex subset $C\subseteq X$
is ball generated, i.\,e. it can be written as an intersection of finite unions of closed balls, formally: there 
exists $\A\subseteq \B$ such that $\bigcap\A=C$, where
\begin{equation*}
\B:=\set*{\bigcup_{i=1}^nB_{r_i}(x_i):n\in \N,\,r_1,\dots,r_n>0,\,x_1,\dots,x_n\in X}.
\end{equation*}
The ball topology $b_X$ is defined to be the coarsest topology on $X$ with respect to which every ball $B_r(x)$ 
is closed. A basis for $b_X$ is given by $\set*{X\sm B:B\in \B}\cup\set*{X}$, where $\B$ is as above. Obviously,
$X$ has BGP if and only if every closed, bounded, convex subset of $X$ is also closed with respect to $b_X$.\par
Ball generated sets and the ball topology were introduced by Godefroy and Kalton in \cite{godefroy} but the notions
implicitly appeared before in \cite{corson}. By \cite{godefroy}*{Theorem 8.1}, every weakly compact subset of a Banach 
space is ball generated. In particular, every reflexive space has the BGP. $c_0$ is an example of a nonreflexive space 
with BGP (see for instance the more general result \cite{basu}*{Theorem 4} on $c_0$-sums). A standard example for a Banach 
space which fails to have the BGP is $\ell^1$ (see the remark at the end of \cite{corson}).\par
We now list some easy remarks on the ball topology (see \cite{godefroy}*{p.197}; some of them may be used later without further notice):
\begin{enumerate}[\upshape(i)]
\item For every $y\in X$, the map $x\mapsto x+y$ is continuous with respect to $b_X$.
\item For every $\lambda>0$, the map $x\mapsto \lambda x$ is continuous with respect to $b_X$.
\item $b_X$ is not a Hausdorff topology, but it is a $T_1$-topology (i.\,e. singletons are closed).
\end{enumerate}
It follows from \cite{godefroy}*{Theorem 8.3} that $X$ has the BGP if and only if the ball topology and the weak topology coincide on $B_X$.
For further information on the ball topology, the BGP and related notions, the reader is referred to \cites{basu,chen,godefroy,granero,linBL}
and references therein.\par
In the paper \cite{basu} by S. Basu many stability results for the BGP are established, in particular, for any family $(X_i)_{i\in I}$ 
of Banach spaces and any $p\in (1,\infty)$, the $\ell^p$-sum $\big[\bigoplus_{i\in I}X_i\big]_p$ has BGP if and only if each $X_i$ has BGP
(\cite{basu}*{Theorem 7}). An analogous result holds for $c_0$-sums (\cite{basu}*{Theorem 4}).\par
In this paper we will study the BGP for direct sums of two spaces only, but with respect to more general norms. We start by recalling the 
necessary definitions: a norm $\norm{\cdot}_E$ on $\R^2$ is called absolute if $\norm{(a,b)}_E=\norm{(\abs{a},\abs{b})}_E$ for all $(a,b)\in \R^2$, 
and it is called normalised if $\norm{(1,0)}_E=\norm{(0,1)}_E=1$. We write $E$ for the normed space $(\R^2,\norm{\cdot}_E)$. For example, the 
standard $p$-norm $\norm{\cdot}_p$ is an absolute, normalised norm for any $p\in [1,\infty]$. Some important properties of absolute, normalised 
norms are listed below (see \cite{bonsall}*{p. 36, Lemma 1 and 2}):
\begin{enumerate}[\upshape(i)]
\item $\norm{(a,b)}_{\infty}\leq\norm{(a,b)}_E\leq\norm{(a,b)}_1 \ \ \forall (a,b)\in \R^2,$
\item $\abs{a}\leq\abs{c},\ \abs{b}\leq\abs{d} \ \Rightarrow \ \norm{(a,b)}_E\leq\norm{(c,d)}_E,$
\item $\abs{a}<\abs{c},\ \abs{b}<\abs{d} \ \Rightarrow \ \norm{(a,b)}_E<\norm{(c,d)}_E.$
\end{enumerate}
For two Banach spaces $X$ and $Y$, their direct sum $X\oplus_E Y$ with respect to $\norm{\cdot}_E$ is defined as the space $X\times Y$
endowed with the norm $\norm{(x,y)}_E:=\norm{(\norm{x},\norm{y})}_E$ for $x\in X$ and $y\in Y$. This is again a Banach space and 
convergence in $X\oplus_E Y$ is equivalent to coordinatewise convergence. For $\norm{\cdot}_E=\norm{\cdot}_p$ one obtains the usual 
$p$-direct sum of Banach spaces.\par
We are going to prove that $X\oplus_E Y$ has the BGP if $X$ and $Y$ have the BGP and the norm $\norm{\cdot}_E$ is G\^ateaux-differentiable 
at $(1,0)$ and $(0,1)$. To do so, we will use a description of absolute, normalised norms by the boundary curve of their unit ball, which
will be discussed in the next section.

\section{Boundary curves of unit balls of absolute norms}\label{sec:boundary}
The following Proposition is quite probably well-known (moreover, its assertion is intuitively clear) but since the author was
not able to find a reference, a formal proof is included here for the readers' convenience.
\begin{proposition}\label{prop:boundary}
Let $\norm{\cdot}_E$ be an absolute, normalised norm on $\R^2$. Then for every $x\in (-1,1)$ there exists exactly one
$y\in (0,1]$ such that $\norm{(x,y)}_E=1$.
\end{proposition}

\begin{Proof}
Let $x\in (-1,1)$. Since the function $t\mapsto \norm{(x,t)}_E$ is continuous with $\lim_{t\to \infty}\norm{(x,t)}_E=\infty$ 
and $\norm{(x,0)}_E=\abs{x}<1$, it follows that there exists $y>0$ such that $\norm{(x,y)}_E=1$. We also have $y\leq \norm{(x,y)}_E=1$.\par
Now  we prove the uniqueness assertion. By symmetry it suffices to consider the case $x\geq 0$. Suppose there exist $0<y_1<y_2\leq 1$ such 
that $\norm{(x,y_1)}_E=\norm{(x,y_2)}_E=1$. Let $0<\lambda<1-y_1/y_2$. It follows that $z:=(x,y_2)+\lambda((1,0)-(x,y_2))=(x+\lambda(1-x),y_2(1-\lambda))$ 
still lies in $B_E$.\par
But $x+\lambda(1-x)>x$ and $y_2(1-\lambda)>y_1$, thus by property (iii) of absolute norms listed in the introduction we must have 
$\norm{z}_E>\norm{(x,y_1)}_E=1$, which is a contradiction.
\end{Proof}

We denote by $f_E$ the function from $(-1,1)$ to $(0,1]$ which assigns to each $x\in (-1,1)$ the corresponding value $y$ given by
Proposition \ref{prop:boundary}. Thus $\norm{(x,f_E(x))}_E=1$ for every $x\in (-1,1)$. The function $f_E$ will be called the upper 
boundary curve of the unit ball $B_E$.\par 
The following properties of $f_E$ are easily verified: $f_E$ is a concave (and hence continuous), even function on $(-1,1)$ with 
$f_E(0)=1$. Further, $f_E$ is increasing on $(-1,0]$ and decreasing on $[0,1)$. In particular, the limits $\lim_{x\nearrow 1}f_E(x)$
and $\lim_{x\searrow -1}f_E(x)$ exist. Thus we may extend $f_E$ to a continuous function from $[-1,1]$ to $[0,1]$, which will be again
denoted by $f_E$.\par
It is possible to characterise properties of the norm $\norm{\cdot}_E$ by corresponding properties of the function $f_E$. As examples
we state below characterisations of strict convexity and strict monotonicity. Once again, this is probably well-known and so the 
(anyway easy) proofs are omitted, but let us first recall the definitions.\par
A Banach space $X$ is strictly convex if $\norm{x+y}=2$ and $\norm{x}=\norm{y}=1$ implies $x=y$.\par
An absolute, normalised norm $\norm{\cdot}_E$ on $\R^2$ is said to be strictly monotone if the following holds:
whenever $a,b,c,d\in \R$ with $\abs{a}\leq\abs{c}$ and $\abs{b}\leq\abs{d}$ and one these inequalities is strict, 
then $\norm{(a,b)}_E<\norm{(c,d)}_E$.
\begin{proposition}\label{prop:convex}
Let $\norm{\cdot}_E$ be an absolute, normalised norm on $\R^2$.\\
The space $E:=(\R^2,\norm{\cdot}_E)$ is strictly convex if and only if $f_E$ is strictly concave\footnote{This means 
$f_E(\lambda x+(1-\lambda)y)>\lambda f_E(x)+(1-\lambda)f_E(y)$ for all $\lambda\in (0,1)$ and all $x,y\in (-1,1)$ with 
$x\neq y$.} on $(-1,1)$ and $f_E(1)=0$.\\
The norm $\norm{\cdot}_E$ is strictly monotone if and only if $f_E$ is strictly decreasing on $[0,1)$ and $f_E(1)=0$.
\end{proposition}
Next we would like to study the smoothness of $\norm{\cdot}_E$ in terms of differentiability of $f_E$. This, too, is quite 
probably known, but the author could not find a reference. Since these results are important for our main result on sums of 
spaces with the BGP, we will provide them here with complete proofs.\par
First recall that, since $f_E$ is concave on $(-1,1)$, it possesses left and right derivatives $f^{\prime}_{E-}$ and $f^{\prime}_{E+}$
on $(-1,1)$ which are decreasing and satisfy $f^{\prime}_{E+}\leq f^{\prime}_{E-}$. Moreover, for every $x_0\in (-1,1)$ and $a\in \R$ we have
\begin{equation}\label{eq:2.1}
f_E(x)\leq f_E(x_0)+a(x-x_0) \ \ \forall x\in (-1,1) \ \Leftrightarrow \ f^{\prime}_{E+}(x_0)\leq a\leq f^{\prime}_{E-}(x_0).
\end{equation}
Also, $f_E$ is differentiable at $x_0\in (-1,1)$ if and only if $f^{\prime}_{E+}$ is continuous at $x_0$ if and only if $f^{\prime}_{E-}$
is continuous at $x_0$. All this follows immediately from the corresponding well-known facts for convex functions, see for example \cite{royden}*{p.113ff.}.\par
For $x\in [-1,1]$, we will denote by $S_E(x)$ the set of support functionals at $(x,f_E(x))$, i.\,e. $S_E(x):=\set*{g\in E^*:\norm{g}_{E^*}=1=g(x,f(x))}$.
\begin{proposition}\label{prop:support}
Let $x_0\in (-1,1)$. For all $a\in [f^{\prime}_{E+}(x_0),f^{\prime}_{E-}(x_0)]$ we have $f_E(x_0)\geq ax_0+1$ and $S_E(x_0)$ consists
exactly of the functionals $g$ of the form
\begin{equation}\label{eq:2.2}
g(x,y)=\frac{ax-y}{ax_0-f_E(x_0)}
\end{equation}
for some $a\in [f^{\prime}_{E+}(x_0),f^{\prime}_{E-}(x_0)]$.
\end{proposition}

\begin{Proof}
Let $a\in [f^{\prime}_{E+}(x_0),f^{\prime}_{E-}(x_0)]$. By \eqref{eq:2.1} we have $f_E(x_0)-ax_0\geq f_E(0)=1$.\par
If $g$ is defined by \eqref{eq:2.2} then it follows from \eqref{eq:2.1} that $g(x,f_E(x))\leq 1$ for all $x\in (-1,1)$.
From this it is easy to deduce that $g(x,y)\leq 1$ for all points $(x,y)$ of norm 1, thus $\norm{g}_{E^*}\leq 1$. Moreover, 
$g(x_0,f_E(x_0))=1$, so $g\in S_E(x_0)$.\par
Conversely, suppose that $g$ is a functional belonging to $S_E(x_0)$. It is of the form $g(x,y)=Ax+By$ for constants $A$ and $B$.
We then have 
\begin{equation}\label{eq:2.3}
Ax\pm Bf_E(x)\leq 1 \ \ \forall x\in (-1,1) \ \ \text{and} \ \ Ax_0+Bf_E(x_0)=1.
\end{equation}
We first prove that $B>0$. If $B\leq 0$, then \eqref{eq:2.3} implies $Ax_0\geq 1$. In the case $x_0>0$ we would obtain, by \eqref{eq:2.3}, 
$1\geq Ax-Bf_E(x)\geq Ax\geq x/x_0$ for all $x\in (0,1)$, which is a contradiction. A similar argument works for $x_0<0$.\par
So we must have $B>0$ and hence it follows from \eqref{eq:2.3} that 
\begin{equation*}
f_E(x)\leq\frac{1}{B}-\frac{A}{B}x \ \ \forall x\in (-1,1).
\end{equation*}
Since $Ax_0+Bf_E(x_0)=1$ we conclude
\begin{equation*}
f_E(x)\leq f_E(x_0)-\frac{A}{B}(x-x_0) \ \ \forall x\in (-1,1).
\end{equation*}
Now \eqref{eq:2.1} implies that $a:=-A/B$ lies in $[f^{\prime}_{E+}(x_0),f^{\prime}_{E-}(x_0)]$.\par 
From $Ax_0+Bf_E(x_0)=1$ we obtain $B=1/(f_E(x_0)-ax_0)$ and hence $A=a/(ax_0-f_E(x_0))$. Thus $g$ is of the form \eqref{eq:2.2}.
\end{Proof}

As is well-known, the norm of a Banach space is G\^{a}teaux-differentiable at a point of norm one if and only if this point has a
unique support functional, which is then the G\^{a}teaux-derivative of the norm at this point. Thus the following is an immediate 
corollary of Proposition \ref{prop:support}.
\begin{corollary}\label{cor:smooth}
Let $\norm{\cdot}_E$ be an absolute, normalised norm on $\R^2$ and $x_0\in (-1,1)$. The norm $\norm{\cdot}_E$ is G\^{a}teaux-differentiable
at $(x_0,f_E(x_0))$ if and only if $f_E$ is differentiable at $x_0$. In this case, the G\^{a}teaux-derivative of $\norm{\cdot}_E$ is given
by
\begin{equation*}
(x,y)\mapsto \frac{f_E^{\prime}(x_0)x-y}{f_E^{\prime}(x_0)x_0-f_E(x_0)}.
\end{equation*}
\end{corollary}
It remains to characterise the support functionals at the end points $(-1,f_E(-1))$ and $(1,f_E(1))$. This requires to distinguish a number of cases.
We will state the result below for completeness, but skip the proof (once again, it should be already known).
\begin{proposition}\label{prop:endpoints}
Let $\norm{\cdot}_E$ be an absolute, normalised norm on $\R^2$. Let $a:=\inf_{x\in [0,1)}f^{\prime}_{E-}(x)\in [-\infty,0]$.\par
For $A,B\in \R$ denote by $g_{A,B}$ the functional given by $g_{A,B}(x,y)=Ax+By$. The following holds:
\begin{enumerate}[\upshape(i)]
\item If $f_E(1)>0$ then $\norm{\cdot}_E$ is G\^{a}teaux-differentiable at each point $(1,b)$ with $b\in (-f_E(1),f_E(1))$ and the 
G\^{a}teaux-derivative at each such point is $g_{1,0}$.
\item $f_E(1)=1$ if and only if $a=0$ if and only if $\norm{\cdot}_E=\norm{\cdot}_{\infty}$. In that case 
$S_E(1)=\set*{g_{A,B}:A,B\geq 0 \ \text{and} \ A+B=1}$.
\item If $a=-\infty$ then $\norm{\cdot}_E$ is G\^{a}teaux-differentiable at $(1,f_E(1))$ with $S_E(1)=\set*{g_{1,0}}$.
\item If $f_E(1)>0$ and $-\infty<a<0$, then $g_{A,B}\in S_E(1)$ if and only if $(A,B)=(\frac{c}{c-f_E(1)},\frac{-1}{c-f_E(1)})$ for 
some $c\in (-\infty,a]$ or $(A,B)=(1,0)$.
\item If $f_E(1)=0$ and $-\infty<a<0$, then $g_{A,B}\in S_E(1)$ if and only if $(A,B)=(1,\pm\frac{1}{c})$ for some $c\in (-\infty,a]$
or $(A,B)=(1,0)$.
\end{enumerate}
\end{proposition}
By symmetry arguments, an analogous characterisation holds for the left endpoint $(-1,f_E(-1))$. Let us also remark that similar characterisations
of support functionals of absolute, normalised norms (on $\C^2$ even) can be found for example in \cite{bonsall}*{p.38, Lemma 4}. These characterisations 
do not use the function $f_E$, but rather the function $\psi$ given by $\psi(t)=\norm{(1-t,t)}_E$ for $t\in [0,1]$.

\section{Direct sums of spaces with the BGP}\label{sec:sums}
We start with the following analogue of \cite{basu}*{Lemma 5}.

\begin{lemma}\label{lemma:balltop}
Let $\norm{\cdot}_E$ be an absolute, normalised norm on $\R^2$ with the following property:
\begin{equation}\label{eq:3.1}
\forall \eps>0 \ \exists s_0>\eps \ \forall s\geq s_0 \ \norm{(1,s-\eps)}_E<s.
\end{equation}
Let $X, Y$ be Banach spaces and $Z:=X\oplus_E Y$. Let $((x_i,y_i))_{i\in I}$ be a net in $B_Z$ which is 
convergent to $0$ in the ball topology $b_Z$. Then $(y_i)_{i\in I}$ converges to $0$ in the topology $b_Y$.\\
Likewise, if $\norm{\cdot}_E$ satisfies
\begin{equation}\label{eq:3.2}
\forall \eps>0 \ \exists s_0>\eps \ \forall s\geq s_0 \ \norm{(s-\eps,1)}_E<s,
\end{equation}
one can conclude that $(x_i)_{i\in I}$ converges to $0$ with respect to $b_X$.
\end{lemma}

\begin{Proof}
The proof is also analogous to that of \cite{basu}*{Lemma 5}. We suppose that $y_i\not\to 0$ with respect to $b_Y$.
Then, by passing to a subnet if necessary, we may assume that there are $y\in Y$ and $r>0$ such that $y_i\in B_r(y)$ 
for all $i\in I$ and $0\in Y\setminus B_r(y)$, i.\,e. $\norm{y}>r$.\par
Put $\eps:=\norm{y}-r$. By \eqref{eq:3.1} we can find $s>\max\set*{\eps,\norm{y}}$ such that $t:=\norm{(1,s-\eps)}_E<s$.\par
Now if $u\in B_X$ and $v\in B_{s-\eps}(sy/\norm{y})$, then by the monotonicity of $\norm{\cdot}_E$,
\begin{equation*}
\norm{(u,v)-(0,sy/\norm{y})}_E=\norm{(\norm{u},\norm{v-sy/\norm{y}})}_E\leq\norm{(1,s-\eps)}_E=t,
\end{equation*}
in other words: $B_X\times B_{s-\eps}(sy/\norm{y})\ssq B_t((0,sy/\norm{y}))$.\par
But for $w\in B_r(y)$ we have
\begin{equation*}
\norm{w-sy/\norm{y}}\leq\norm{w-y}+\norm{y-sy/\norm{y}}\leq r+s-\norm{y}=s-\eps,
\end{equation*}
thus $B_r(y)\ssq B_{s-\eps}(sy/\norm{y})$.\par
Altogether it follows that $(x_i,y_i)\in B_t((0,sy/\norm{y}))$ for every $i\in I$. But $0\not\in B_t((0,sy/\norm{y}))$,
since $t<s$. So the complement of $B_t((0,sy/\norm{y}))$ is a $b_Z$-neighbourhood of $0$ not containing any of the points
$(x_i,y_i)$. With this contradiction the proof is finished.
\end{Proof}

As mentioned in the introduction, $X$ has the BGP if and only if the ball topolgy and the weak topology of $X$ coincide 
on $B_X$ (\cite{godefroy}*{Theorem 8.3}). Thus we can, as in \cite{basu}, derive the following stability result.
\begin{corollary}\label{cor:sums}
Let $\norm{\cdot}_E$ be an absolute, normalised norm on $\R^2$ satisfying both \eqref{eq:3.1} and \eqref{eq:3.2}.
Let $X$ and $Y$ be Banach spaces with the BGP. Then $X\oplus_E Y$ also has the BGP.
\end{corollary}

\begin{Proof}
It follows from Lemma \ref{lemma:balltop} that for every bounded net $((x_i,y_i))_{i\in I}$ in $X\oplus_E Y$
which is convergent to some point $(x,y)$ in the ball topology we also have $x_i\to x$ and $y_i\to y$ in the 
resprective ball topologies of $X$ and $Y$. Since $X$ and $Y$ have the BGP, it follows that these nets also 
converge in the weak topology of $X$ resp. $Y$, which in turn implies $(x_i,y_i)\to (x,y)$ in the weak topology 
of $X\oplus_E Y$. Thus $X\oplus_E Y$ has the BGP.
\end{Proof}

It remains to determine which absolute norms satisfy the conditions \eqref{eq:3.1} and \eqref{eq:3.2}. As it turns
out, \eqref{eq:3.1} resp. \eqref{eq:3.2} is equivalent to the G\^{a}teaux-differentiablility of $\norm{\cdot}_E$ at
$(0,1)$ resp. $(1,0)$. To prove this we will use the description of the norm by its upper boundary curve $f_E$ from 
the previous section and the following version of the mean value theorem for one-sided derivatives (see for instance 
\cite{saksbook}*{p.204} or \cite{walter}*{p.358} for an even more general statement).
\begin{theorem}\label{thm:meanvalue}
Let $I$ be an interval and $f:I \rightarrow \R$ a continuous function. Let $J$ be another interval. Suppose that the 
right derivative $f_{+}^{\prime}(x)$ exists and lies in $J$ for all but at most countably many interior points from $I$.
Then
\begin{equation*}
\frac{f(b)-f(a)}{b-a}\in J \ \ \forall a,b\in I \ \text{with} \ a\neq b.
\end{equation*}
An analogous statement holds for the left derivative.
\end{theorem}

\begin{proposition}\label{prop:gateaux}
Let $\norm{\cdot}_E$ be an absolute, normalised norm on $\R^2$. $\norm{\cdot}_E$ is G\^{a}teaux-differentiable at $(0,1)$
resp. $(1,0)$ if and only if \eqref{eq:3.1} resp. \eqref{eq:3.2} holds.
\end{proposition}

\begin{Proof}
We only prove the statement for $(0,1)$, the other case follows from this one by considering instead of $\norm{\cdot}_E$ the 
norm given by $\norm{(x,y)}_F:=\norm{(y,x)}_E$.\par
Assume first that $\norm{\cdot}_E$ is G\^{a}teaux-differentiable at $(0,1)$. By Corollary \ref{cor:smooth} the function $f_E$ 
is differentiable at $0$ and the G\^{a}teaux-derivative of $\norm{\cdot}_E$ at $(0,1)$ is given by
\begin{equation*}
(x,y)\mapsto -f_E^{\prime}(0)x+y.
\end{equation*}
But this G\^{a}teaux-derivative must be the projection onto the second coordinate, thus $f_E^{\prime}(0)=0$.\par
For each real number $s>0$ we define $f_s(x):=sf_E(x/s)$ for $x\in (-s,s)$. The functions $f_s$ are continuous
and differentiable from the right with $f_{s+}^{\prime}(x)=f_{E+}^{\prime}(x/s)$.\par
Let $\eps>0$. Since $f_{E+}^{\prime}$ is continuous at 0 (cf. the remarks preceding Proposition \ref{prop:support}) we 
can find $\delta\in (0,1)$ such that $\abs*{f_{E+}^{\prime}(x)}<\eps$ for every $x\in (-\delta,\delta)$. Let 
$s_0>\max\set*{\eps,1/\delta}$ and $s\geq s_0$. Then $\abs*{f_{s+}^{\prime}(x)}<\eps$ for all $x\in (0,1)$ and thus 
by Theorem \ref{thm:meanvalue} $\abs*{f_s(1)-f_s(0)}<\eps$, hence $f_s(1)>s-\eps$.\par
This implies $\norm{(1,s-\eps)}_E<s$, for otherwise we would have $s=\norm{(1,f_s(1))}_E\geq\norm{(1,s-\eps)}_E\geq s$, so
$\norm{(1,s-\eps)}_E=s$, which would mean $f_E(1/s)=1-\eps/s$ and thus we would obtain the contradiction $f_s(1)=s-\eps$.
This completes one direction of the proof.\par
To prove the converse we assume that \eqref{eq:3.1} holds but $\norm{\cdot}_E$ is not G\^{a}teaux-differentiable at $(0,1)$.
Then by Corollary \ref{cor:smooth}, the function $f_E$ is not differentiable at 0. Since $f_E$ is increasing on $(-1,0]$ we 
have $a:=f_{E-}^{\prime}(0)\geq 0$ and because $f_E$ is even we have $f_{E+}^{\prime}(0)=-a$. Hence $a>0$ and by \eqref{eq:2.1}
$f_E(x)\leq f_E(0)+f_{E+}^{\prime}(0)x=1-ax$ for all $x\in (-1,1)$.\par
If we define $f_s$ as above it follows that
\begin{equation}\label{eq:3.3}
f_s(x)\leq s-ax \ \ \forall x\in (-s,s), \forall s>0.
\end{equation}
By \eqref{eq:3.1} we can choose $s>\max\set*{1,a}$ such that $\norm{(1,s-a)}_E<s$. Then by \eqref{eq:3.3} $f_s(1)\leq s-a$
and hence $s=\norm{(1,f_s(1))}_E\leq \norm{(1,s-a)}_E<s$. This contradiction finishes the proof.
\end{Proof}

Putting Corollary \ref{cor:sums} and Proposition \ref{prop:gateaux} together we obtain the final result.
\begin{corollary}
Let $\norm{\cdot}_E$ be an absolute, normalised norm on $\R^2$ which is G\^{a}teaux-differentiable at $(0,1)$ and $(1,0)$.
Let $X$ and $Y$ be Banach spaces with the BGP. Then $X\oplus_E Y$ also has the BGP.
\end{corollary}
This result contains in particular the case of $p$-sums for $1<p\leq\infty$ that---as we mentioned in the introduction---was 
already treated in \cite{basu} (even for infinite sums). As was also mentioned in \cite{basu}, the BGP cannot be 
stable under infinite $\ell^1$-sums (since $\ell^1$ itself does not have the BGP), but it is open whether $X\oplus_1 Y$ has the 
BGP whenever $X$ and $Y$ have it.

\address
\email

\end{document}